\documentclass[11pt]{article}

\usepackage{a4wide}
\usepackage[english]{babel}
\usepackage{amssymb}
\usepackage{amsthm}
\usepackage{enumerate}
\usepackage{amsmath}

\newtheorem{thm}{Theorem}

\newtheorem{proposition}[thm]{Proposition}

\newcommand{\N}{\mathbb N}
\newcommand{\Z}{\mathbb Z}
\newcommand{\R}{\mathbb R}

\renewcommand{\L}{\mathbb L}

\newcommand{\ind}{\mathchoice{\mathrm {1\mskip-4.1mu l}} {\mathrm{ 1\mskip-4.1mu l}} {\mathrm {1\mskip-4.6mu l}} {\mathrm {1 \mskip-5.2mu l}}}

\newcommand{\F}{{\mathcal F}}
\newcommand{\A}{{\mathcal A}}
\newcommand{\B}{{\mathcal B}}
\newcommand{\C}{{\mathcal C}}
\newcommand{\esp}{\mathbb E}

\newcommand{\cqfd}{\begin{flushright}$\square$\end{flushright}}

\begin{document}

 \title{\huge\textbf{Independence of Four Projective Criteria for\\
 the Weak Invariance Principle}}
 \author{Olivier Durieu\footnote{E-mail: olivier.durieu@etu.univ-rouen.fr.\newline
web page: http://www.univ-rouen.fr/LMRS/Persopage/Durieu/\newline
Laboratoire de Mathématiques 
Rapha\"el Salem,
UMR 6085 CNRS-Universit\'e de Rouen.}}

\maketitle

\begin{abstract}
Let $(X_i)_{i\in\Z}$ be a regular stationary process for a given filtration.
The weak invariance principle holds under the condition $\sum_{i\in\Z}\|P_0(X_i)\|_2<\infty$ 
(see \cite{Han79}, \cite{DedMer03}, \cite{DedMerVol07}). In this paper, we show that this criterion is independent of other known criteria:
the martingale-coboundary decomposition of Gordin (see \cite{Gor69}, \cite{Gor73}), the criterion of Dedecker and Rio (see \cite{DedRio00})
and the condition of Maxwell and Woodroofe (see \cite{MaxWoo00}, \cite{PelUte05}, \cite{Vol06}, \cite{Vol07}).


\end{abstract}

\section{Introduction and Results}

The aim of this paper is to study the relation between several criteria to get the weak invariance principle of Donsker
in dependent case.
In the paper \cite{DurVol07}, the independence between three of them is already shown. These criteria are the
martingale-coboundary decomposition, the projective criterion of Dedecker and Rio and the Maxwell-Woodroofe condition.
Here, we consider a fourth one ($\sum_{i\in\Z}\|P_0(X_i)\|_2<\infty$) and we show that it is independent of the three others.
Let us begin by the statements of the four criteria.

\medskip

Let $(\Omega,\A,\mu)$ be a probability space and let $T$ be a bijective bimeasurable transformation of $\Omega$ preserving $\mu$. We assume $(\Omega,\A,\mu, T)$ is an ergodic dynamical system. Let $f:\Omega\longrightarrow\R$ be
a measurable function with zero mean. We recall that the process $(f\circ T^i)_{i\in\N}$ satisfies the weak invariance principle if the process
$$
\frac{1}{\sqrt{n}}\sum_{k=0}^{\lfloor tn\rfloor-1}f\circ T^k,\,t\in[0,1]
$$
converges in distribution 
to a Gaussian process in the space $D([0,1])$ provided with the Skorohod topology (see Billingsley \cite{Bil68}). 
Let $\F$ be a sub-$\sigma$-algebra of $\A$ such that $T^{-1}\F\subset \F$. We denote by $\F_i$ the $\sigma$-algebra $T^{-i}\F$. 
The function $f$ is called regular with respect to the filtration $(\F_i)_{i\in\Z}$ if
$$
\esp(f|\F_{-\infty})=0\mbox{ and }\esp(f|\F_{+\infty})=f.
$$
In the sequel, we assume that $f$ is a square integrable function and we write $\L^p$ for $\L^p(\mu)$, $p\geq 1$.

\medskip

$\bullet$ The first criterion is the martingale-coboundary decomposition due to Gordin \cite{Gor69}.
We will restrict our attention to the martingale-coboundary decomposition in $\L^1$ (see \cite{Gor73}, and 
\cite{EssJan85} for a complete proof).
We say that $f$ admits such a decomposition if $f=m+g-g\circ T$
where $(m\circ T^i)_{i\in\Z}\subset \L^1$ is a martingale difference sequence and $g\in\L^1$.
If $m\in\L^2$, then the central limit theorem holds. Further, if 
$\frac{1}{\sqrt{n}}\max_{i\leq n}|g\circ T^i|$
goes to $0$ in probability, the weak invariance principle holds (see Hall and Heyde \cite{HalHey80}).
If $f$ is a regular function with respect to the filtration $(\F_i)_{i\in\Z}$ 
then the martingale-coboundary decomposition in $\L^1$ is equivalent to
\begin{equation}\label{1}
\sum_{i=0}^\infty\esp(f\circ T^i|\F_0)\mbox{ and }\sum_{i=0}^\infty f\circ T^{-i}-\esp(f\circ T^{-i}|\F_0)\mbox{ converge in }\L^1,
\end{equation}
see Voln\'y \cite{Vol93}.
Remark that if the process $(f\circ T^i)_{i\in\Z}$
is adapted to $(\F_i)_{i\in\Z}$, the second sum is equal to zero.

\medskip

$\bullet$ The Dedecker and Rio criterion is satisfied if
\begin{equation}\label{2}
\sum_{k=1}^\infty f\esp(f\circ T^k|\F_0) \mbox{ converges in }\L^1.
\end{equation}
According to Dedecker and Rio \cite{DedRio00}, in the adapted case, this condition implies the weak invariance principle.

\medskip

$\bullet$ The Maxwell-Woodroofe condition (see \cite{MaxWoo00}) is satisfied if 
\begin{equation}\label{3}
\sum_{n=1}^\infty\frac{\|\esp(S_n(f)|\F_0)\|_2}{n^{\frac{3}{2}}}<+\infty
\end{equation}
where $S_n(f)=\sum_{i=0}^{n-1}f\circ T^i$.
In the adapted case, Peligrad and Utev \cite{PelUte05} proved that this condition implies the weak invariance principle. In the general case, the weak invariance principle
 holds as soon as (\ref{3}) and
\begin{equation*}
\sum_{n=1}^\infty\frac{\|S_n(f)-\esp(S_n(f)|\F_n)\|_2}{n^{\frac{3}{2}}}<+\infty,
\end{equation*}
(see Voln\'y \cite{Vol06}, \cite{Vol07}).
\medskip

The independence between these three criteria is proved in \cite{DurVol07}. Here we add a new criterion.
Let us denote by $H_k=\L^2(\F_k)$ the space of $\F_k$-measurable functions which are square integrable
and denote by $P_k$ the orthogonal projection operator onto the space $H_k\ominus H_{k-1}$. 
For $f\in\L^2$,
$$
P_k(f)=\esp(f|\F_k)-\esp(f|\F_{k-1}).
$$

\medskip

$\bullet$ Let $f$ be a regular function for the filtration $(\F_i)_{i\in\Z}$.
As a consequence of a result given by Heyde \cite{Hey74} (see \cite{Vol93} Theorem 6)
the central limit theorem holds as soon as 
\begin{equation}\label{4}
\sum_{i\in\Z}\|P_0(f\circ T^i)\|_2<\infty.
\end{equation}
In the adapted case, this result and the weak invariance principle were proved by Hannan \cite{Han73}, \cite{Han79} under the assumption that $T$ is weakly mixing. 
Hannan's weak invariance principle was proved without the extra assumption by Dedecker and Merlev\`ede \cite{DedMer03}, Corollary 3. Finally, in the general case, the weak invariance principle under 
(\ref{4}) is due to Dedecker, Merlev\`ede and Voln\'y \cite{DedMerVol07}, Corollary 2.

\bigskip

Our main result is the following theorem.

\begin{thm}\label{th}
Conditions \textnormal{(\ref{1})}, \textnormal{(\ref{2})}, \textnormal{(\ref{3})} and
\textnormal{(\ref{4})} are pairwise independent:
in all ergodic dynamical system with positive entropy,
for each couple of conditions among the four, there exists an $\L^2$-function  satisfying the first condition
but not the second one.
\end{thm}

With the results of \cite{DurVol07}, it remains to prove the independence of (\ref{1}), (\ref{2}), (\ref{3}) with (\ref{4}).

\section{Proof of Theorem \ref{th}}

Let $(\Omega,\A,\mu, T)$ be an ergodic dynamical system with entropy greater or equal than $1$.
Let $\B$ and $\C$ be two independent sub-$\sigma$-algebra of $\A$. Let $(e_i)_{i\in\Z}$ be a sequence of independent
$\B$-measurable random variables in
 $\{-1,1\}$ such that $\mu(e_i=-1)=\mu(e_i=1)=\frac{1}{2}$
and  $e_i=e_0\circ T^i$, $i\in\Z$. We denote by $\F_0$ the $\sigma$-algebra generated by $\C$ and $e_i$ for $i\leq 0$ and we set
$\F_i=T^{-i}\F_0$.
Note that the case of entropy in $(0,1)$ can be studied by using another Bernoulli shift.

We introduce three sequences with the following properties:

$(\theta_k)_{k\in\N}\subset (0,+\infty)$;
 
$(\rho_k)_{k\in\N}\subset(0,1)$ such that $\sum_{k\geq 0}\rho_k<1$;
 
$(N_k)_{k\in\N}\subset\N$ such that $N_{k+1}> N_k$.

We can always find a sequence $(\varepsilon_k)_{k\in\N}\subset(0,1)$ such that $\sum_{k\geq0}\theta_kN_k\sqrt{\varepsilon_k}<\infty$. 
So, we fix such a sequence and
denote by $f$ the function defined by
$$
f=\sum_{k=1}^\infty\theta_ke_{-N_k}\ind_{A_k}
$$
where the sequence of sets $(A_k)_{k\in\N}$ verifies:
\begin{itemize}
\item $\forall k\in \N$, $A_k\in\C$,
\item the sets $A_k$ are disjoint,
\item $\exists a\in(0,1)$, $\forall k\in\N$, $a\rho_k\leq\mu(A_k)\leq\rho_k$,
\item $\forall k\in\N$, $\forall i\in\{0,\dots ,N_k\}$, $\mu(T^{-i}A_k\Delta A_k)\leq\varepsilon_k$.
\end{itemize}
The construction of these sets is done in detail in \cite{DurVol07}.
 
\medskip 
 
First, remark that the process $(f\circ T^i)_{i\in\Z}$ is adapted to
the filtration $(\F_i)_{i\in\Z}$ (then $f$ is regular) and
$$
f\in \L^2 \mbox{ if and only if } \sum_{k\geq 1}\theta_k^2\rho_k<\infty.
$$
The next proposition is proved in \cite{DurVol07}.
\begin{proposition}\label{prop}
For the function $f$ previously defined,
\begin{enumerate}[\rm i.]
\item 
$
\textnormal{(\ref{1})}\quad\Leftrightarrow\quad\sum_{k\geq 1}\theta_k\rho_k\sqrt{N_k}<\infty;
$
\item 
$
\textnormal{(\ref{2})}\quad\Leftrightarrow\quad\sum_{k\geq 1}\theta_k^2\rho_k\sqrt{N_k}<\infty;
$
\item 
$
\textnormal{(\ref{3})}\quad\Leftrightarrow\quad
\sum_{n\geq 1} n^{-\frac{3}{2}}\left(\sum_{k\geq 1}\theta_k^2\min(n,N_k)\rho_k\right)^{\frac{1}{2}}<\infty.
$
\end{enumerate}
\end{proposition}

We can state an analogous result for condition (\ref{4}).

\begin{proposition}\label{sp}
For the function $f$ previously defined,
$$
\sum_{i\in \Z}\|P_0(f\circ T^i)\|_2<\infty\mbox{ if and only if }\sum_{k\geq 1}\theta_k\sqrt{\rho_k}<\infty.
$$
\end{proposition}
This proposition is proved in Section \ref{preuve}.

\subsection*{Counterexamples}

Now we give two counterexamples proving Theorem \ref{th}.

\begin{enumerate}
\item We consider the function $f$ defined by the sequences $\theta_k=\frac{2^k}{k}$,
$\rho_k=\frac{1}{4^k}$, $N_k=k$.
Then $f\in\L^2$ and using Proposition \ref{prop}, we get:
\begin{enumerate}[a)]
\item $\sum_{k\geq 1}\theta_k\rho_k\sqrt{N_k}=\sum_{k\geq 1}\frac{1}{2^k\sqrt{k}}<\infty$
and then (\ref{1}) is verified.
\item $\sum_{k\geq 1}\theta_k^2\rho_k\sqrt{N_k}=\sum_{k\geq 1}k^{-\frac{3}{2}}<\infty$
and then (\ref{2}) is verified.
\item $\sum_{k\geq 1}\theta_k^2\min(n,N_k)\rho_k
=\sum_{k\geq 1}\frac{\min(n,k)}{k^2}$.
But 
$\sum_{k=1}^n\frac{1}{k}\leq1+\ln(n)$ and $\sum_{k=n+1}^\infty\frac{n}{k^2}\leq 1$.
Then
$
\sum_{n\geq 1} n^{-\frac{3}{2}}\left(\sum_{k\geq 1}\theta_k^2\min(n,N_k)\rho_k\right)^{\frac{1}{2}}
\leq \sum_{n\geq 1} n^{-\frac{3}{2}}\sqrt{\ln(n)+2}<\infty
$
and (\ref{3}) is verified.
\item $\sum_{k\geq 1}\theta_k\sqrt{\rho_k}=\sum_{k\geq 1}\frac{1}{k}$ diverges and then by Proposition \ref{sp}, 
(\ref{4}) is not satisfied.
\end{enumerate}

This counterexample shows that  none of the conditions (\ref{1}), (\ref{2}) and (\ref{3}) implies (\ref{4}).

\medskip

\item Now, if we consider the function $f$ defined by the sequences: 
$\theta_k=\frac{2^k}{k^{\frac{3}{2}}}$, $\rho_k=\frac{1}{4^k}$, $N_k=2^{4k}$, then $f\in\L^2$ and:
\begin{enumerate}[a)]
\item $\sum_{k\geq 1}\theta_k\rho_k\sqrt{N_k}=\sum_{k\geq 1}\frac{2^k}{k^{\frac{3}{2}}}$ diverges and then
(\ref{1}) does not hold.
\item $\sum_{k\geq 1}\theta_k^2\rho_k\sqrt{N_k}=\sum_{k\geq 1}\frac{2^{2k}}{k^3}$ diverges and then
(\ref{2}) does not hold.
\item $\sum_{k\geq 1}\theta_k^2\min(n,N_k)\rho_k=\sum_{k\geq 1}\frac{\min(n,2^{4k})}{k^3}
\geq n\sum_{k\geq\lfloor\frac{\ln n}{4\ln 2}\rfloor}\frac{1}{k^3}\geq 8\ln^2 2\frac{n}{\ln^2 n}$.
Then 
$
\sum_{n\geq 1} n^{-\frac{3}{2}}\left(\sum_{k\geq 1}\theta_k^2\min(n,N_k)\rho_k\right)^{\frac{1}{2}}
\geq 8\ln^2 2\sum_{n\geq 1}\frac{1}{n\ln n}
$ 
diverges and (\ref{3}) does not hold.
\item On the other hand, $\sum_{k\geq 1}\theta_k\sqrt{\rho_k}=\sum_{k\geq 1}\frac{1}{k^\frac{3}{2}}<\infty$ and then 
(\ref{4}) holds.
\end{enumerate}

This shows that (\ref{4}) does not imply any conditions (\ref{1}), (\ref{2}) or (\ref{3}).
\end{enumerate}

\subsection*{Remarks}
\begin{enumerate}[i.]
\item In fact, we showed a little more than Theorem \ref{th}. We got that the three conditions  
(\ref{1}), (\ref{2}) and (\ref{3}) together are independent of (\ref{4}).
\item To show that (\ref{4}) does not imply (\ref{1}), it is enough to consider a linear process
$f=\sum_{i\in\Z}a_i\xi_i$ where $(\xi_i)_{i\in\Z}$ is an iid sequence with $\mu(\xi_0=1)=\mu(\xi_0=-1)=\frac{1}{2}$
and $a_i=\frac{1}{i^{\frac{3}{2}}}$.
\end{enumerate}

\section{Proof of Proposition \ref{sp}}\label{preuve}

First of all, $(f\circ T^i)_{i\in\Z}$ is adapted to the filtration and then for all $i<0$,
$$
P_0(f\circ T^{i})=0.
$$
For $i\geq 0$, we have
\begin{eqnarray}
P_0(f\circ T^i)&=&\esp(f\circ T^i|\F_0)-\esp(f\circ T^i|\F_{-1})\nonumber\\
&=&\sum_{k\geq 1}\theta_k\left[\esp(e_{-N_k+i}|\F_0)-\esp(e_{-N_k+i}|\F_{-1})\right]\ind_{A_k}\circ T^i.\nonumber
\end{eqnarray}
Since $e_j$ is $\F_0$-measurable for $j\leq 0$ and independent of $\F_0$ for $j>0$,
\begin{eqnarray}
\esp(e_{-N_k+i}|\F_0)-\esp(e_{-N_k+i}|\F_{-1})&=&
\left\{\begin{array}{ll}
e_0&{\rm if}\quad i=N_k\\
0&{\rm otherwise}
\end{array}\right.\nonumber\\
&=&e_0\ind_{\{i=N_k\}}.\nonumber
\end{eqnarray}
Thus
\begin{eqnarray}
P_0(f\circ T^i)&=&\sum_{k\geq 1}\theta_ke_0\ind_{\{i=N_k\}}\ind_{A_k}
+\sum_{k\geq 1}\theta_ke_0\ind_{\{i=N_k\}}(\ind_{T^{-i}A_k}-\ind_{A_k})\nonumber\\
&=&\qquad I_1(i)\qquad+\qquad I_2(i)\qquad.\nonumber
\end{eqnarray}

For $I_2$, we use the fact that $\mu(A_k\Delta T^{-i}A_k)\leq \varepsilon_k$ for $0\leq i\leq N_k$ to get
\begin{eqnarray}
\|I_2(i)\|_2&\leq&\sum_{k\geq 1}\theta_k\ind_{\{i=N_k\}}\|e_0\ind_{A_k\Delta T^{-i}A_k}\|_2\nonumber\\
&\leq&\sum_{k\geq 1}\theta_k\ind_{\{i=N_k\}}\sqrt{\varepsilon_k}.\nonumber
\end{eqnarray}
Remark for each $i\geq 0$, there is at most one integer $k$ such that $N_k=i$ and
for each $k\geq 1$,
there exists an integer $i$ such that $i=N_k$. 
We deduce
\begin{eqnarray}
\sum_{i\geq 0}\|I_2(i)\|_2&\leq&
\sum_{i\geq 0}\sum_{k\geq 1}\theta_k\ind_{\{i=N_k\}}\sqrt{\varepsilon_k}\nonumber\\
&=&\sum_{k\geq 1}\theta_k\sqrt{\varepsilon_k}\nonumber
\end{eqnarray}
which is finite by the assumptions.

\medskip

Thus, $\sum_{i\geq 0}\|P_i(f)\|_2$ is converging if and only if $\sum_{i\geq 0}\|I_1(i)\|_2$ is converging.
Now for a fixed $i$, since the sets $A_k$ are disjoint and since there is at most one $k$ such that $N_k=i$,
we have
\begin{eqnarray}
\|I_1(i)\|_2&=&\sqrt{\sum_{k\geq 1}\theta_k^2\ind_{\{i=N_k\}}\mu(A_k)}\nonumber\\
&=&\sum_{k\geq 1}\theta_k\ind_{\{i=N_k\}}\sqrt{\mu(A_k)}.\nonumber
\end{eqnarray}
Finally,
\begin{eqnarray}
\sum_{i\geq 0}\|I_1(i)\|_2&=&
\sum_{i\geq 0}\sum_{k\geq 1}\theta_k\ind_{\{i=N_k\}}\sqrt{\mu(A_k)}\nonumber\\
&=&\sum_{k\geq 1}\theta_k\sqrt{\mu(A_k)}.\nonumber
\end{eqnarray}
We can conclude the proof using  $a\rho_k\leq \mu(A_k)\leq \rho_k$.
\cqfd

\subsection*{Acknowledgment}
I would like to thank J\'er\^ome Dedecker for suggesting this problem.

\small
\bibliographystyle{plain}

\end{document}